\documentclass[12pt,a4paper]{article} 
\usepackage{amsfonts}
\newtheorem{proposition}{Proposition} 
\begin{document}
\title{A construction of representations for quantum groups: 
an example of ${\cal U}_q(\mathfrak{so}(5))$}
\author{Pavel \v S\v tov\'\i\v cek\\
\small Department of Mathematics\\
\small Faculty of Nuclear Science, CTU\\
\small Trojanova 13, 120 00 Prague, Czech Republic\\
\small stovicek@kmdec.fjfi.cvut.cz
}
\date{}
 \def\CC{{\cal C}}
 \def\UU{{\cal U}}
 \def\FF{{\cal F}}
 \def\MM{{\cal M}}
 \def\RR{{\cal R}}
 \def\BC{\mathbb{ C}}
\def\BZ{\mathbb{ Z}}
 \def\vphi{\varphi}
 \def\eps{\epsilon}
 \def\veps{\varepsilon}
 \def\sg{\sigma}
 \def\gtg{\mathfrak{ g}}
 \def\gk{\mathfrak{ k}}
 \def\gh{\mathfrak{h}}
 \def\gX{\mathfrak{ X}}
 \def\gso{\mathfrak{ so}}
 \def\gR{\mathfrak{ R}}
\newcommand{\be}{\begin{equation}}
\newcommand{\ee}{\end{equation}}
\newcommand{\bea}{\begin{eqnarray}}
\newcommand{\eea}{\end{eqnarray}}
 \def\lo{_{(1)}}
 \def\lt{_{(2)}}
 \def\loo{_{(1)(1)}}
 \def\lot{_{(1)(2)}}
 \def\lto{_{(2)(1)}}
 \def\ltt{_{(2)(2)}}
 \def\Zp{\BZ_{+}}
 \def\Uq{{\cal U}_q}
 \def\Lp{L^{+}}
 \def\Lm{L^{-}}
 \def\Lpm{L^{\pm}}
 \def\Lmp{L^{\mp}}
 \def\Z{Z^{\ast}} 
 \def\z{z^{\ast}} 
 \def\qhI{q^{H_1}}
 \def\qhII{q^{H_2}} 
 \def\xpI{X^+_1}
 \def\xpII{X^+_2}
 \def\xmI{X^-_1} 
 \def\xmII{X^-_2} 
 \def\tx{\tilde\xi}
 \def\tp{\tilde\varphi}
 \def\te{\tilde\varepsilon}
 \def\td{\tilde\Delta}
 \def\sgn{\mbox{sgn}}
 \def\spn{\mbox{span}}
 \def\tr{\mbox{tr}}
 \def\diag{\mbox{diag}}
 \def\id{\mbox{id}}
 \def\Funq{\mbox{Fun}_{q}}
 \def\idr{\langle{\cal R}\rangle} 
 \def\half{{1\over2}}
\maketitle

\begin{abstract}
A short description is given of a construction of representations for quantum 
groups. The method uses infinitesimal dressing transformation on quantum 
homogeneous spaces and is illustrated on an example of 
$\UU_q(\gso(5))$. 
\end{abstract}

\section{Introduction}

The purpose of this paper is to illustrate a construction of representations on an 
explicit example, namely the deformed enveloping algebra $\UU_q(\gso(5))$. 
We are going to describe the construction as well, however, its detailed 
presentation will appear elsewhere. The basic ingredient is the infinitesimal 
dressing transformation 
on a quantum homogeneous space, in  analogy with the celebrated method of 
orbits due to Kirillov and Kostant. 

The construction generalizes and simplifies some results derived in the papers 
\cite{JMP1,Grass, Coher, JAlg,JMP2} and also \cite{Bulg,JSch}. Let us mention just 
a few additional papers dealing also with constructions of representations 
of quantum groups and/or 
with quantum homogeneous spaces \cite{PW,TT,BL,DDF}, but taking a different 
point of view or applying other methods. 

Concerning the deformation parameter, we assume that $q>0$, $q\neq1$. All 
fractional powers of $q$ are supposed to be positive. 

\section{Construction} 

We assume that we are given a
bialgebra $\UU$ with the counit denoted by $\veps$ and the
comultiplication denoted by $\Delta$,  and a unital algebra
$\CC$. Moreover, $\CC$ is supposed to be a left $\UU$-module with the
action denoted by $\xi$, and fulfilling two conditions:
\bea
&& \xi_x\cdot 1=\veps(x)\,1,\quad \forall x\in\UU, \\
&&  \xi_x\cdot(fg)=(\xi_{x\lo}\cdot f)(\xi_{x\lt}\cdot g),\quad
\forall x\in\UU,\ \forall f,g\in\CC.
\eea
If convenient we shall write $\xi(x)\cdot f$ instead of
$\xi_x\cdot f$. The second condition (2) is nothing but Leibniz
rule. Here and everywhere in what follows we use Sweedler's notation:
$\Delta x=x\lo\otimes x\lt$. 

\begin{proposition}
Suppose that a linear mapping $\vphi:\UU\to\CC$ satisfies $\vphi(1)=1$ and
\be
\vphi(xy)=\bigl(\xi_{x\lo}\cdot\vphi(y)\bigr)\vphi(x\lt),
\quad\forall x,y\in\UU.
\ee 
Then the prescription
\be
x\cdot f:=(\xi_{x\lo}\cdot f)\,\vphi(x\lt),
\quad\forall x\in\UU,\ \forall f\in\CC,
\ee 
defines a left $\UU$-module structure on $\CC$ and it holds
\be
x\cdot(fg)=(\xi_{x\lo}\cdot f)(x\lt\cdot g),\quad
\forall x\in\UU,\ \forall f,g\in\CC.
\ee 
Particularly,
\be
\vphi(x)=x\cdot 1,\quad\forall x\in\UU.
\ee
\indent Conversely, suppose that
$\UU\otimes\CC\to\CC:x\otimes f\mapsto x\cdot f$
is a left $\UU$-module structure on $\CC$ such that the rule (5) is
satisfied. Then the linear mapping $\vphi:\UU\to\CC$ defined by the
equality (6) fulfills (3), and consequently the prescription (4) holds true. 
\end{proposition}

Let us suppose, as usual, that $\UU$
is generated as an algebra by a set of generators $\MM\subset\UU$. Let
$\FF$ be the free algebra generated by $\MM$. Thus $\UU$ is identified
with a quotient $\FF/\idr$ where $\idr$ is the ideal generated by a
set of defining relations $\RR\subset\FF$. Let $\pi$ be the factor
morphism, $\pi:\FF\to\UU$.
We set $\te:=\veps\circ\pi$ and
\be
\tx_x\cdot f:=\xi_{\pi(x)}\cdot f,\quad\forall x\in\FF,\
\forall f\in\CC .
\ee
In addition we impose the following condition on the set of generators $\MM\subset\UU$: 
\be 
\Delta(\MM)\subset\spn_{\BC}(\MM_1\otimes\MM_1)\quad
\mbox{where}\quad\MM_1:=\MM\cup\{1\}.
\ee 
Then it is natural to define a comultiplication $\td$ on $\FF$ by the equality
$\td(x_1\dots x_n):=\Delta(x_1)\dots\Delta(x_n)$, $x_i\in\MM$.
As $\UU$ is a bialgebra $\idr$ must be, at the same time, a coideal.

It is not difficult to check that $\FF$ becomes this way a bialgebra and
that the triple $(\FF,\tx,\CC)$ fulfills the original conditions (1)
and (2), just replacing $\UU$ with $\FF$ and $\xi$ with $\tx$.
One finds that to any mapping 
$\vphi:\MM\to\CC$ there exists a unique linear extension
$\tp:\FF\to\CC$ such that $\tp(1)=1$  and the property
\be 
\tp(xy)=\bigl(\tx_{x\lo}\cdot\tp(y)\bigr)\tp(x\lt),
\ee 
is satisfied for all $x,y\in\FF$.

The final step in the construction is to decide when the mapping
$\tp$ can be factorized from $\FF$ to $\UU=\FF/\idr$.
\begin{proposition} 
Suppose that there is given a mapping $\vphi:\MM\to\CC$ and let
$\tp$ be its extension to $\FF$ as described above. If
\be 
(\pi\otimes\tp)\circ\td(\RR)=0
\ee 
then $\tp(\idr)=0$ and so there exists a unique linear mapping
$\vphi':\UU\to\CC$ such that $\tp=\vphi'\circ\pi$. Moreover, $\vphi'=1$ 
and $\vphi'$ satisfies the condition (3).\newline  
\indent The same conclusions hold true provided $\RR$ fulfills a 
stronger condition than that of being a coideal, namely 
\be 
\td(\RR)\subset\idr\otimes\FF+\FF\otimes\FF\RR, 
\ee 
and $\tp$ satisfies a weaker condition 
\be 
\tp(\RR)=0.
\ee 
\end{proposition}

Particularly this construction goes through for the standard deformed enveloping 
algebras $\UU=\Uq(\gtg)$ in the FRT description \cite{FRT} where $\gtg$ 
is any simple complex Lie algebra from the four principal series 
$A_\ell$, $B_\ell$, $C_\ell$ and $D_\ell$. So the generators are arranged
in respectively upper and lower triangular matrices $L^+$ and $L^-$, and 
the set $\RR$ is given by the usual RLL relations. 

On the other hand the unital algebra $\CC$ is generated by quantum 
antiholomorphic coordinate functions $\z_{jk}$, $j<k$, on the generic dressing 
orbit of dimension $(\dim_{\BC}\gtg-\mbox{rank}\,\gtg)/2$. The elements are 
arranged in an upper triangular matrix $Z$ with units on the diagonal, and the 
defining relations are given in terms of its Hermitian adjoint $\Z$, namely 
\be 
R_{12}\Z_2Q\Z_1Q^{-1} = \Z_1Q\Z_2Q^{-1}R_{12} 
\ee
where $Q$ is the diagonal part of the R-matrix $R$. 

The infinitesimal dressing transformation $\xi$ is prescribed on the generators, 
\be 
\xi(\Lp_1)\cdot \Z_2 = R_{21}^{-1}\Z_2Q, \quad
\xi(\Lm_1)\cdot\Z_2 = \Z_1Q\Z_2Q^{-1}(\Z_1)^{-1}.
\ee
It can be extended to an arbitrary element from $\CC$ with the aid of 
Leibniz rule (2).
The mapping $\vphi$ is defined on the generators as well,
\be 
\vphi(\Lp)=D^{-1},\quad
\vphi(\Lm)=\Z D^2(\Z)^{-1}D^{-1}
\ee 
where $D$ is an arbitrary complex diagonal matrix obeying the conditions
\be 
\det(D)=1\quad\mbox{and}\quad K_{12}D_1D_2=K_{12}. 
\ee 
Here $K$ is a matrix related to the R-matrix via the equality
\be 
R_{12}-R_{21}^{-1}=(q-q^{-1})(P-K_{12}), 
\ee 
$P$ stands for the flip operator. 

\section{Example: $\Uq(\gso(5))$}

We shall use the Drinfeld--Jimbo description of  $\Uq(\gso(5))$ \cite{D,J}, 
with the six generators $\qhI,\ \qhII,\ \xpI,\ \xpII,\ \xmI,\ \xmII$, the relations 
\bea
&  [\,q^{H_1} ,\, q^{H_2}\,] = 0,\nonumber\\ 
&  q^{H_1}\,{X_1}^\pm = q^{\pm 1}\,q^{H_1}\,{X_1}^\pm,\ 
  q^{H_1}\,{X_2}^\pm = q^{\mp 1}\,q^ {H_1}\,{X_2}^\pm,\nonumber\\ 
&  q^{H_2}\,{X_1}^\pm = q^{\mp 1}\,q^{H_2}\,{X_1}^\pm,\
  q^{H_2}\,{X_2}^\pm = q^{\pm 2}\, q^{H_2}\,{X_2}^\pm,\nonumber\\ 
&  [\,{X_1}^+,\,{X_1}^-\,] = \frac{q^{H_1} - q^{-H_1}}{q-q^{-1}},\
  [\,{X_2}^+,\,{X_2}^-\,] = \frac{q^{H_2} - q^{-H_2}}{q-q^{-1}},  \\
&  [\,\xpI ,\, \xmII\,] = 0,\  [\,\xpII ,\, \xmI\,] = 0,\nonumber\\ 
&  ({X_2}^\pm)^2\, {X_1}^\pm - 
    \left( q^{-1} + q \right) \,{X_2}^\pm\,{X_1}^\pm\,{X_2}^\pm 
    +{X_1}^\pm\,({X_2}^\pm)^2\,= 0,\nonumber\\ 
& \begin{array}{ll}
({X_1}^\pm)^3\,{X_2}^\pm - 
    \left(  q^{-1}+1+ q \right) \,({X_1}^\pm)^2\,{X_2}^\pm\,{X_1}^\pm\, 
      \\
    \mbox{      }
 + \left( q^{-1}+1+ q \right) \, {X_1}^\pm\,{X_2}^\pm\, ({X_1}^\pm)^2
     - {X_2}^\pm\,({X_1}^\pm)^3\,= 0 \,.
\end{array} \nonumber
\eea
and the comultiplication 
\be
\Delta(q^{H_i})= q^{H_i}\otimes q^{H_i},\ 
\Delta(X_i^\pm)=X_i^\pm\otimes q^{-{1\over 2}\,H_i}+
q^{{1\over 2}\,H_i}\otimes  X_i^\pm,\quad i=1,2. 
\ee
One can pass from the FRT description to the Drinfeld--Jimbo generators using the 
equalities 
\bea
&  \Lp_{11}=\Lm_{55}=q^{H_1+H_2},\  \Lp_{22}=\Lm_{44}=q^{H_1},\nonumber\\ 
& \Lp_{12}=(q-q^{-1})\,q^{-1/2}\,\xmII\, q^{H_1+\half H_2}, \ 
 \Lp_{23}=(q-q^{-1})\,q^{-1/2}\,\xmI\, q^{\half H_1},\\ 
 & \Lm_{34}=(q-q^{-1})\,q^{1/2}\,\xpI\, q^{\half H_1}, \ 
 \Lm_{45}=(q-q^{-1})\,q^{1/2}\,\xpII\, q^{H_1+\half H_2} . \nonumber
\eea

Only 4 among the 10 generators $\z_{jk}$, $1\le j<k\le 5$, are independent. We denote 
the independent generators by $w_1,\dots,w_4$ and make the following choice:  
\be 
\z_{12}=w_1,\ \z_{13}=w_2,\ \z_{14}=w_3,\ \z_{23}=w_4.
\ee 
The remaining entries can be expressed in terms of $w_1,\dots,w_4$ as well,
\bea
& \z_{45}= - w_1,\ \z_{34}= -q^{1/2}\, w_4,\ 
\z_{35}= - q^{-1/2}\, w_2+q^{1/2}\,w_1 w_4 ,\nonumber\\ 
& \z_{24}=-\frac{q^{1/2}}{1+q} w_4^{\: 2},\ 
 \z_{15}= - w_1 w_3 - \frac{q^{-1/2}}{1+q} w_2^{\: 2},\\ 
& \z_{25}=-q^{-1}\, w_3 - q^{-1/2}\, w_2 w_4 + 
  \frac{q^{1/2}}{1+q} w_1 w_4^{\: 2} \,.\nonumber
\eea
The algebra $\CC$ is then determined by the relations 
\bea
&  w_2 w_1=q\,w_1 w_2,\  w_3 w_2=q\,w_2 w_3,\  
 w_3 w_4=q\,w_4 w_3,\\ \nonumber 
&  w_3 w_1=w_1 w_3 -q^{-1/2}\,(q-1)\,w_2^{\: 2},\ 
 w_4 w_2=w_2 w_4 -q^{-1/2}\,(q-q^{-1})\,w_3,\\  
& w_4 w_1=q^{-1}\, w_1 w_4 + (1-q^{-2})\,w_2. \nonumber
 \eea   
Consequently, the ordered monomials 
$w_1^{\: n_1}w_2^{\: n_2}w_3^{\: n_3}w_4^{\: n_4}$, $n_1,n_2,n_3,n_4\in\Zp$,   
form an algebraic basis of $\CC$. 

The infinitesimal dressing transformation is prescribed on the generators as follows: 
\bea 
&& \xi(\qhI)\cdot \{w_1,w_2,w_3,w_4\}=\{q^{-1}\,w_1,\, w_2,\, q\,w_3,\, q\,w_4\},
\nonumber\\ 
&& \xi(\qhII)\cdot \{w_1,w_2,w_3,w_4\}=\{q^2\,w_1,\, q\,w_2,\, w_3,\, q^{-1}\,w_4\},
\nonumber\\ 
&& \xi(\xmI)\cdot \{w_1,w_2,w_3,w_4\}=
  \{0,\,-q^{1/2}\,w_2,\, q^{1/2}\, w_3,\, -1\}, \\ 
&& \xi(\xmII)\cdot \{w_1,w_2,w_3,w_4\}=\{-q^{1/2},\,0,\,0,\,0\}, \nonumber \\ 
&& \xi(\xpI)\cdot \{w_1,w_2,w_3,w_4\}=
 \left\{-q^{-1/2}\,w_2,\, q^{-1/2}\,w_3,\, 0,\, q^{1/2}\,w_3,\, 
 \frac{q^{1/2}}{1+q}\,w_4^{\:2}\right\},  \nonumber\\ 
&& \xi(\xpII)\cdot \{w_1,w_2,w_3,w_4\} \nonumber\\
&& \qquad\qquad\qquad\quad 
 =\left\{q^{-1/2}\,w_1^{\:2},\, w_1 w_2,\,  -\frac{1}{1+q}\,w_2^{\:2},\, 
  -q^{-1}\,w_1 w_4+q^{-2}\,w_2 \right\}.\nonumber
\eea 

Let us turn to the mapping $\vphi$. The constraints (16) imply that 
\be
D=\diag\left(q^{\half\,\sg_1+\sg_2},q^{\half\,\sg_1},1,
  q^{-\half\,\sg_1},q^{-\half\,\sg_1-\sg_2}\right)  
\ee 
where $\sg_1,\sg_2\in\BC$ are parameters. 
A straightforward calculation gives 
\bea
& \vphi(\qhI)=q^{-\half\,\sg_1},\ \vphi(\qhII)=q^{-\sg_2},\ 
\vphi(\xmI)=\vphi(\xmII)=0,\nonumber\\ 
& \vphi(\xpI)=-\frac{q^{\half-{1\over 4}\,\sg_1}}{1+q}\,[\sg_1]_{q^{1/2}}\, w_4,\ 
\vphi(\xpII)=-q^{-\half(1+\sg_2)}\,[\sg_2]_q\, w_1, 
\eea
where $[m]_p:=(p^m-p^{-m})/(p-p^{-1})$. 

The final step is to calculate the modified action according to the prescription (4). 
Here is the result:
\bea
&& \qhI\cdot w_1^{\:n_1}w_2^{\:n_2}w_3^{\:n_3}w_4^{\:n_4} =
  q^{-n_1+n_3+n_4-\half\,\sg_1}\, w_1^{\:n_1}w_2^{\:n_2}w_3^{\:n_3}w_4^{\:n_4}, 
 \nonumber\\ 
&& \qhII\cdot w_1^{\:n_1}w_2^{\:n_2}w_3^{\:n_3}w_4^{\:n_4} =
  q^{2\,n_1+n_2-n_4-\sg_2}\, w_1^{\:n_1}w_2^{\:n_2}w_3^{\:n_3}w_4^{\:n_4}, 
 \nonumber\\ 
&& \xmI\cdot w_1^{\:n_1}w_2^{\:n_2}w_3^{\:n_3}w_4^{\:n_4}  = \nonumber\\
&& \qquad\qquad
        -\, q^{\half(-n_1+n_2-n_3-n_4)+{1\over4}\,\sg_1}\, [n_2]_{q^{1/2}}
  \, w_1^{\:n_1+1}w_2^{\:n_2-1}w_3^{\:n_3}w_4^{\:n_4} \nonumber\\
&& \qquad\qquad  
       +\, q^{\half(-n_1+n_3-n_4)+{1\over4}\,\sg_1}\, [n_3]_q
  \, w_1^{\:n_1}w_2^{\:n_2+1}w_3^{\:n_3-1}w_4^{\:n_4} \nonumber\\
&& \qquad \qquad
         -\, q^{\half(-n_1+n_3)+{1\over4}\,\sg_1}\, [n_4]_{q^{1/2}}
  \, w_1^{\:n_1}w_2^{\:n_2}w_3^{\:n_3}w_4^{\:n_4-1}, \nonumber\\
&& \xmII\cdot w_1^{\:n_1}w_2^{\:n_2}w_3^{\:n_3}w_4^{\:n_4}  = 
 \nonumber\\
&& \qquad\qquad 
        -\, q^{\half(1-n_2+n_4+\sg_2)}\, [n_1]_q  
  \, w_1^{\:n_1-1}w_2^{\:n_2}w_3^{\:n_3}w_4^{\:n_4}, \nonumber\\
&& \xpI\cdot w_1^{\:n_1}w_2^{\:n_2}w_3^{\:n_3}w_4^{\:n_4}  = \\
&& \qquad\qquad
        -\, q^{-1+\half(n_1-n_3-n_4)+{1\over4}\,\sg_1}\, [n_1]_q 
  \, w_1^{\:n_1-1}w_2^{\:n_2+1}w_3^{\:n_3}w_4^{\:n_4} \nonumber\\
&& \qquad\qquad  
       +\, q^{-1+\half(-n_1+n_2-n_3-n_4)+{1\over4}\,\sg_1}\, [n_2]_{q^{1/2}}
  \, w_1^{\:n_1}w_2^{\:n_2-1}w_3^{\:n_3+1}w_4^{\:n_4} \nonumber\\
&& \qquad \qquad
         +\, \frac{q^{\half(1-n_1+n_3)-{1\over4}\,\sg_1}}{1+q}\, 
         [n_4-\sg_1]_{q^{1/2}}
  \, w_1^{\:n_1}w_2^{\:n_2}w_3^{\:n_3}w_4^{\:n_4+1}, \nonumber\\
&& \xpII\cdot w_1^{\:n_1}w_2^{\:n_2}w_3^{\:n_3}w_4^{\:n_4}  = 
\nonumber\\
&& \qquad\qquad
        q^{-\half(1-n_2+n_4+\sg_2)}\, [n_1+n_2-n_4-\sg_2]_q 
  \, w_1^{\:n_1+1}w_2^{\:n_2}w_3^{\:n_3}w_4^{\:n_4} \nonumber\\
&& \qquad \qquad
         -\, \frac{q^{-1+n_1+\half\,n_2+n_3 - {3\over2}(n_4+\sg_2)}}{1+q}\, 
         [n_3]_q 
  \, w_1^{\:n_1}w_2^{\:n_2+2}w_3^{\:n_3-1}w_4^{\:n_4} \nonumber\\
  && \qquad\qquad  
       +\, q^{-1+n_1+\half\,n_2+n_3 -n_4- {3\over2}\,\sg_2}\, [n_4]_{q^{1/2}}
  \, w_1^{\:n_1}w_2^{\:n_2+1}w_3^{\:n_3}w_4^{\:n_4-1} \nonumber\\
&& \qquad \qquad
         -\,(q-1)\, q^{-{5\over2}+n_1+\half(n_2 -n_4)- {3\over2}\,\sg_2}\, 
         [n_4]_{q^{1/2}} [n_4-1]_{q^{1/2}} \nonumber\\
&& \qquad\qquad\qquad   \times\:       
  \, w_1^{\:n_1}w_2^{\:n_2}w_3^{\:n_3+1}w_4^{\:n_4-2}. \nonumber
\eea

Note that $1\in\CC$ is a lowest weight vector ($\xmI\cdot1=\xmII\cdot1=0$), with 
the lowest weight determined by 
$\qhI\cdot1=q^{-\half\,\sg_1}$, $\qhII\cdot1=q^{-\sg_2}$. Consequently, the cyclic 
submodule $\UU\cdot1$ is finite-dimensional and irreducible provided 
$\sg_1,\sg_2\in\Zp$, and this way one can obtain, in principle, all 
finite-dimensional irreducible representations of $\Uq(\gso(5))$. For example, 
if $\sg_1=1$, $\sg_2=0$, then  $\UU\cdot1$ is a 4-dimensional vector space 
spanned by the vectors: 1, $w_4$, $w_2-q\,w_1 w_4$, 
$(1+q)w_3+q^{3/2}\,w_2 w_4$.

\vskip 10pt
 \noindent{\bf Acknowledgements.} Partial support from the grant
 202/96/0218 of Czech Grant Agency is gratefully acknowledged.

\end{document}